\documentclass[12pt]{article}
\usepackage{graphicx}

\newcommand{\C}{{\bf C}}

\newcommand{\R}{{\bf R}}

\newcommand{\Z}{{\bf Z}}

\begin{document}

\begin{large}

\begin{center}

\textbf{{\Large On the \textit{genera} of moment-angle manifolds associated to dual-neighborly polytopes, combinatorial formulas and sequences.}}

\bigskip

Santiago L\'opez de Medrano.

\bigskip

\end{center}

\end{large}
\bigskip

\begin{small}
Abstract. For a family of polytopes of even dimension $2p$, known as \textit{dual-neighborly}, it has been shown for $p\ne 2$ that the associated intersection of quadrics is a connected sum of sphere products $S^p\times S^p$. In this article we give formulas for the number of terms in that connected sum. Certain combinatorial operations produce new polytopes whose associated intersections are also connected sums of sphere products and we give also formulas for their number. These include a large amount of simple polytopes, including many odd-dimensional ones.
\end {small}
\section*{Introduction}
To every simple polytope $P$ there is associated a manifold $Z(P)$ of the same dimension known in different works as its \textit{(real) moment-angle manifold}, \textit{universal abelian cover} ([D-J]), \textit{polyhedral product} ([B-B-C-G]) or \textit{intersection of quadrics (more precisely, of coaxial ellipsoids)} ([LdM], [Go-LdM2]).\\

The topology of $Z(P)$ cannot be described in full generality, but it has been described for some large families of polytopes $P$. One of them is the family of \emph{dual-neighborly} polytopes $P$ of even dimension $2p\ne 4$ for which it was shown in [Gi-LdM] that they are connected sums of copies of the sphere product $S^p\times S^p$. For $P$ of odd dimension $2p+1\ge 5$  it was proved under a certain additional hypothesis (probably unnecessary) that they are connected sums of copies of the sphere product $S^p\times S^{p+1}$. Together with them, for each such polytope $P$ of dimension at least 5 there is an infinite lattice of polytopes obtained by applying iteratively in all possible ways a well-known operation $P\mapsto P'$, known as the \textit{book}\footnote{This name has been used for many years in the theory of intersection of quadrics. In the literature on polytopes this construction is called the \textit{wedge} on $P$.} or the \textit{wedge} construction for which it was shown that the associated manifold is a connected sum of sphere products $S^a \times S^b$ with factors of different dimensions. The number of combinatorially different dual-neighborly polytopes of even dimension grows enormously fast with their dimension and number of facets ([P]) and all their corresponding infinite lattices are disjoint, so we know that a huge part of the simple polytopes have associated manifolds which are connected sums of sphere products.\\

However, only when $p=1$ or when the number of facets of the polytope is at most $2p+3$ do we know the exact number of summands and precisely which sphere products appear. For $p=1$ and $P$ the $n$-gon the number of those summands (i.e., the genus of the surface $Z(P)$) is known to be $2^{n-3}(n-4)+1$. The same sequence of numbers appears in many other geometric and combinatorial questions, see [Sl] and [Go-LdM]. It appeared for the first time (to our knowledge) in a 1935 paper by Coxeter as the genera of surfaces obtained by a certain construction of his ([Co]) and was found independently around 1980 by Hirzebruch (unpublished, but see [Hi]) as the genera of a certain family of real surfaces that are intersections of quadrics. Only much later Coxeter's construction was recognized as a precursor of what is now called a \textit{polyhedral product} and that intersections of quadrics of the type that Hirzebruch had considered are another instance of the same construction. I have seen no evidence that they were ever conscious of that coincidence.\\

We will give now a generalization of this formula for all even dimensions that gives the number of terms in the connected sum, which is natural to call the \emph{genus} of $Z(P)$. The formula is actually valid homologically  even in those cases where it can be conjectured (but not yet proved) that $Z$ is a connected sum of sphere products and can be extended to all the polytopes obtained from $P$ through the book construction. It is still not understood how the combinatorics of $P$ determines the precise products that appear in the connected sum after several applications of the book construction and not only their number.\\

The genus formula follows, in the case of a dual-neighborly polytope of dimension $d=2p$ and $n$ facets, from known combinatorial formulas for the number of the faces of a neighborly polytope in each dimension . One can obtain from them the Euler characteristic of $Z(P)$ and therefore its genus.\\

This direct result is useful for computations, but very messy and not too useful. The real work consisted in the search of a better formula. After several failed attempts, two elements opened the road for a solution. First, the appearance of the sequences of genera for a few small values of $p$ in the \textit{Sloane Encyclopedia of Sequences} ([Sl]), which included generating functions for them that suggested immediately a nice and simple conjecture for all dimensions. Secondly, a specific direct formula for the number of faces of a dual-neighborly polytope in the book by Br\o ndsted ([Br]) that was more suited for our computations. From these facts, a long computation led us to a proof of the conjecture (Theorems 1 and 2).\\

This gave a collateral proof of various combinatorial identities that we had tried to prove in our first attempts. Other by-products are a new interpretation and some new formulas for the cases appearing in the \textit{Sloane Encyclopedia of Sequences}, as well as an infinite generalization of them. One could search for extensions of the various interpretations of those few cases.\\  

Additionally, it was shown in [Gi-LdM] that other geometric operations on the polytopes, such as truncation of vertices, induce in the associated manifolds operations that preserve connected sums of sphere products. We give the genus also for all $P$ obtained from the ones above by iterated vertex truncations and book constructions in any order.\\

So this work is a quantitative continuation of [Gi-LdM], giving the explicit number of sphere products in the connected sums. Alas, this time Sam Gitler was no more among us to participate and enjoy this extension of our work. 

\section{Background.}

The construction of $Z(P)$ for a simple $d$-polytope $P$ with $n$ facets can be described as follows: one can assume that $P$ is embedded in a $d$-dimensional affine subspace $A$ of $\R^{n}$ in such a way that $A\cap \R_{+}^{n}=P$ and $A$ intersects transversely every coordinate subspace of $\R_{+}^{n}$. Then $Z(P)$ is the union of all the images of $P$ under all compositions of reflections of $\R^{n}$ on its coordinate subspaces $\{x_{i}=0\}$. $Z(P)$ is a combinatorial $d$-manifold that can be easily smoothed as a transversal intersection of ellipsoids ([LdM], [Go-LdM2]). $Z(P)$ can also be constructed abstractly as a quotient of $P\times \Z_2^n$ under the identifications in the facets of $P$ corresponding to the fixed points under the reflections on the coordinate subspaces ([LdM],[D-J]).

The book construction consists in taking the product $P\times [0,1]$ and one of its facets $F$ and, for each given point $x\in F$, identifying all points $(x,t)$ for $t\in [0,1]$ into a single point. Under this operation, the dimension and the number of facets of the polytope increases by $1$ and we denote by $P'$ the polytope so obtained and by $Z'$ the corresponding manifold. A geometric construction of $Z'$ and manipulations with homology exact sequences shows that the total homology (i.e., the direct sum of the homology groups) of $Z$ and $Z'$ is the same.\\

One can consider compositions of an arbitrary number of book constructions along different facets. Following [B-B-C-G2] we denote by $P^J$, $Z^J$, where $J=(j_1, j_2, \dots, j_n)$, the result of applying $j_i$ times the book construction on the $i$-th facet of $P$ for $i=1,\dots, n$. See [B-B-C-G2] for details of a more general construction and a combinatorial description of $P^J$ in the dual context of simplicial complexes.\\

The simple polytope $P$ is called \emph{dual-neighborly} if every collection of $k$ facets of $P$ has a non-empty intersection, for all $k\le d/2$ (Cf. [Br], p.92 and [B-M], p.114). They are dual to the much studied \emph{neighborly} ones. It can be proved that $Z(P)$ is  $[d/2-1]$-connected if, and only if, $P$ is dual-neighborly. So, for a $2p$-dimensional (respectively, $(2p+1)$-dimensional dual-neighborly $P$), $Z(P)$ has homology only in dimension $p$ (respectively, in dimensions $p$ and $p+1$), other than the $0$ and top dimensional ones.\\ 

It is known that if two dual-neighborly polytopes have the same dimension $d$ and the same number of facets $n$, then they have the same number of $k$-dimensional faces for all $k$ from $0$ to $d$. Various explicit formulas are known for this number of $k-$faces as a function of $(d,n)$, one of which will be more suited to our purposes. We will give this formula in the next section.\\

Now suppose that $P$ is dual-neighborly polytope and of even dimension $d=2p$ and $n$ facets. Then the homology of $Z(P)$ is free and is non-trivial only in the middle dimension $p$, so it has the homology of a connected sum of copies of the sphere product $S^{p} \times S^{p}$. It was shown in [Gi-LdM] that it is actually diffeomorphic to such a connected sum if $d > 4$, but the number of those products was not given. The book construction applied on any facet of $P$ gives also a dual-neighborly polytope $P'$ of dimension $2p+1$ and the manifold $Z(P')$ is a connected sum of copies of $S^{d} \times S^{d+1}$ (even in the case $d=4$). Since the book construction preserves connected sums of sphere products, after any number of further applications of it one obtains again manifolds that are connected sums of products $S^{a} \times S^{b}$ for various pairs of dimensions $(a,b)$. \\

\section{The Euler characteristic $\chi(Z(P))$ for dual-neighborly polytopes $P$ of even dimension.}

Let $P$ be a dual-neighborly polytope of even dimension $d=2p$ with $n$ facets. Since $n\ge d+1$ (with equality only for the simplex) it is better to use the parameters $p=d/2$ and $m$ defined as $$m=n-d-1=n-2p-1$$that starts with $m=0$.\\

The numbers $f_{k}$ of $k$-faces of $P$ (for $k=0,\dots,d$) are determined only by the numbers $p,m$. Explicit formulas for them can be deduced from the formulas for the number of $k$-faces of a neighborly $d$-polytope with $n$ vertices in [Gr], section 9.2. A direct explicit formula for $f_{k}$ is given in [Br], p.113, which in our notation becomes:

\begin{center}
$ f_{k}=\Sigma_{j=0}^{p}$     ${j} \choose{k}$      ${m + j} \choose{j}$+$\Sigma_{j=0}^{p-1}$ ${2p-j} \choose{k}$                                                         ${m + j} \choose{j}$
\end{center}

Now, the polytope $P$ is reflected on the coordinate hyperplanes of $\R^n$ to give a cell decomposition of $Z$ formed by $2^n$ cells of dimension $d$ which are all copies of $P$. A face of dimension $k$ has $d-k$ coordinates equal to zero so it is reflected only on $n-d+k$ hyperplanes and therefore produces $2^{n-d+k}=2^{m+k+1}$ faces. Thus the total number of $k$-faces of $Z$ is $f_{k}\times 2^{m+k+1}$.\\

Therefore, the Euler characteristic of $Z(P)$, which we denote by  $\chi(p,m))$, is the alternating sum

\begin{center}
$\chi(p,m)=\Sigma_{k=0}^{2p} (-1)^k 2^{m+k+1}\large{(}\Sigma_{j=0}^{p}$     ${j} \choose{k}$      ${m + j} \choose{j}$+$\Sigma_{j=0}^{p-1}$ ${2p-j} \choose{k}$                                                        ${m + j} \choose{j}$ $)$ \\
\end{center}

This formula is useful for computations, even for $d,m$ in the thousands, since it can be easily programmed in the computer. It also shows that, for any fixed $p$, $\chi(p,m)$ is of the form $2^{m+1}$  times a polynomial in $m$ of degree $p$. But otherwise it is quite messy and opaque. For example, it is easy to see directly that $\chi(p,0)=2$ ($P$ is the simplex $\Delta^{2p}$ and $Z(P)$ is the sphere $S^{2p}$), $\chi(p,1)=2(1+(-1)^p)$ ($P$ is $\Delta^{p}\times \Delta^{p}$ and $Z(P)=S^{p}\times S^{p}$) and it is known that $\chi(p,2)=2+(-1)^p(4p+6)$ ([LdM]). But these facts are not clear from the formula.\\

In any case, it is convenient to simplify it: factoring $2^{m+1}$, our formula can be re-arranged as follows:
\begin{center}
$\chi(p,m)=2^{m+1}(\Sigma_{j=0}^{p} \Sigma_{k=0}^{2p} (-2)^k$  ${j} \choose{k}$ ${m + j} \choose{j}$                $+\Sigma_{j=0}^{p-1} \Sigma_{k=0}^{2p} (-2)^k$  ${2p-j} \choose{k}$ ${m + j} \choose{j}$ $)$ \\
\end{center}

Now, since ${j} \choose{k}$ $=0$ if $k>j$, we have:

\begin{center}
$\Sigma_{k=0}^{2p} (-2)^k$  ${j} \choose{k}$  $= \Sigma_{k=0}^{j} (-2)^k$  ${j} \choose{k}$  $=(-2+1)^j =(-1)^j$
\end{center}

And, since ${2p-j} \choose{k}$ $=0$ if $k>2p-j$, we have:

\begin{center}

$\Sigma_{k=0}^{2p} (-2)^k$  ${2p-j} \choose{k}$  $= \Sigma_{k=0}^{2p-j} (-2)^k$  ${2p-j} \choose{k}$  $=(-2+1)^{2p-j} =(-1)^{2p-j}=(-1)^j$
\end{center}

And therefore we obtain a better formula

\begin{center}
$\chi(p,m)=2^{m+1}(\Sigma_{j=0}^{p} (-1)^j$  ${m + j} \choose{j}$                $+\Sigma_{j=0}^{p-1} (-1)^j$ ${m + j} \choose{j}$ $)$ \\
\end{center}

An even simpler formula can be obtained by computing the generating function of the above expression parametrized by $m$ for a fixed $p$:
$$\Sigma_{m\ge 0} \chi(p,m) z^m$$

Since the formula for $\chi(p,m)$ involves two sums which differ only in their length,  we can cover both cases with the general sequence
\begin{center}
$S(r,m)=2 \Sigma_{j=0}^{r} (-1)^j 2^m$ ${m + j} \choose{j}$\\
\end{center}
and its corresponding generating function
$$\Sigma_{m\ge 0} S(r,m) z^m$$
which is the sum for $j=0,\dots r$, of the generating functions
\begin{center}
$2 \Sigma_{m\ge 0} (-1)^j 2^m$ ${m + j} \choose{j}$ $z^m$ 
$=2 (-1)^j \Sigma_{m\ge 0}$ ${m + j} \choose{j}$ $(2z)^m$.\\
\end{center}
Now, it is well known (and easy to prove) that
\begin{center}
$\Sigma_{m\ge 0} $ ${m + j} \choose{j}$ $y^m= \frac{1}   {(1-y)^{j+1}}$\,\,\,\,\,\,\,\,\,\,(*)\\
\end{center}
\noindent which gives
$$\Sigma_{m\ge 0} S(r,m) z^m
=2 \Sigma_{j=0}^{r}(-1)^j \frac{1}   {(1-2z)^{j+1}}=-2\Sigma_{j=0}^{r}  \frac{1}   {(2z-1)^{j+1}}$$

This is a geometric progression with sum

$$-2\frac{     \frac{1}{(2z-1)^{r+2}}-\frac{1}{2z-1}}{\frac{1}{2z-1}-1}=-2\frac{     \frac{1}{(2z-1)^{r+1}}-1}{1-(2z-1)}=-2\frac{     \frac{1}{(2z-1)^{r+1}}-1}{2-2z}$$
$$=\frac{1}   {(1-z)} + \frac{1}   {(z-1)(2z-1)^{r+1}}$$

The generating function for $ \chi(p,m)$ is the sum of two instances of the above expression evaluated a $r=p$ and $r=p-1$,
% Error $$-2 ( \frac{1}   {2z-1} - \frac{1}   {(2z-2)(2z-1)^{p+1}}+ \frac{1}   {2z-1} - \frac{1}   {(2z-2)(2z-1)^{p}})$$
which add up to
$$ \frac{2}  {1-z} + \frac{1} {(z-1)}  \left(\frac{1}{(2z-1)^{p+1}}+\frac{1}{(2z-1)^p}\right)$$
$$= \frac{2}   {1-z} + \frac{2z}   {(z-1)(2z-1)^{p+1}}$$

We can obtain another expression for the general term of this series, by working each part separately:\\

The first term of the function $ \frac{2}   {1-z}$ is $ \Sigma_{m\ge 0} 2 z^m$.\\

And the second term has two factors: $2z/(z-1)$ and $1/(2z-1)^{p+1}$. The first one is simply
$ \frac {2z}   {z-1} $ so in its series the coefficient of $z^i$ is $-2$ for $i>1$
and for the second factor (using formula (*) again):
\begin{center}
$ \frac{1}   {(2z-1)^{p+1}}= \frac{(-1)^{p+1}}   {(1-2z)^{p+1}}=
 (-1)^{p+1} \Sigma_{j\ge 0} $ ${j + p} \choose{p}$ $(2z)^j$
\end{center}

so the coefficient of $z^j$ is   $(-1)^{p+1}$ ${j + p} \choose{p}$ $2^j$\\

So, in the product, the coefficient of $z^m$ is
\begin{center}
$\Sigma_{j=0}^{m-1}  (-1)^p 2^{j+1}$ ${j + p} \choose{p}$
\end{center}
since for $j\ge m$ the coefficient of the first factor is $0$.\\

To which we still have to add the first term. So, finally, the coefficient of $z^m$, which is $ \chi(p,m)$, is\\
\begin{center}
$ \chi(p,m)= 2+ (-1)^p \Sigma_{j=0}^{m-1}2^{j+1}$ ${j + p} \choose{p}$.\\
\end{center}
We have proved:\\

\textbf{Theorem 1.} \textit{If $P$ is a dual-neighborly polytope of even dimension $d=2p$ and $n=d+m+1$ facets, then $\chi(p,m)$, the Euler characteristic of $Z(P)$, can be expressed in any of the following equivalent forms:}\\

(i) $\chi(p,m)=2^{m+1}(\Sigma_{j=0}^{p} (-1)^j$  ${m + j} \choose{j}$     $+\Sigma_{j=0}^{p-1} (-1)^j$ ${m + j} \choose{j}$ $)$\\

(ii) $ \chi(p,m)=(-1)^p \Sigma_{j=0}^{m-1} 2^{j+1}$ ${j + p} \choose{p}$ $+2$\\

(iii) \textit{$\chi(p,m)$, as a sequence parametrized by $m$, has generating function}

$$ \frac {2} {1-z} + \frac {2z} {(z-1)(2z-1)^{p+1}}$$

These formulas are valid for any $p$ and $m$, and in the cases where we know that $Z(P)$ is a connected sum of sphere products $S^p\times S^p$, we can derive the number of them from the formulas.\\

But these formulas do not extend to the manifolds obtained by the book construction on $P$: just the first application gives an odd dimensional manifold with $\chi=0$. And also, the even dimensional ones obtained by iteration of the book construction on $P$ may include products of two odd-dimensional spheres that contribute with negative terms to $\chi$, so $\chi$ will not depend only on the dimension and number of facets of the corresponding polytope.\\

We will solve these problems by the introduction of the concept of \textit{genus} of such a connected sum.  

\section{The \textit{genus} $g(Z(P))$ for dual-neighborly polytopes $P$ of even dimension and associated polytopes.}

For a connected sum of sphere products $M$, we can naturally define its \textit{genus} $g(M)$, as the number of products in the sum, as in the case of surfaces.\\

Let $\beta(M)$ the sum of the Betti numbers $\beta_i(M)$, then $\beta(M)=2 g(M)+2$, or
$$g(M)= \frac{\beta(M)}{2}-1$$

This definition can be extended to any manifold or even any topological space with finite $\beta$ by the above formula. In some cases this may not be an integer (if $X$ is a point then $g(X)=-1/2$). Some properties of this generalized genus (for example, that it is additive with respect to the connected sum of manifolds) will be considered elsewhere.\\

In the case that all the summands are of the form $S^p \times S^p$,  there is a simple relation between $g(M)$ and the Euler characteristic of $M$: 
$$\chi(M)=2+ (-1)^p 2 g(M)$$
or
$$g(p,m)= (-1)^p  \chi(p,m)/2 -(-1)^p$$
and again, this relation is valid for any manifold with the same homology groups. This includes the manifolds $Z(P)$ where $P$ is a dual-neighborly polytope and for which the formulas for the Euler characteristic are still valid. So with this extension of the concept of genus we can state:\\

\textbf{Theorem 2.} \textit{If $P$ is a dual-neighborly polytope of even dimension $d=2p$ and $n=d+m+1$ facets, then $g(p,m)$, the genus of $Z(P)$, can be expressed in any of the following equivalent forms:}\\

(i) $g(p,m)=(-1)^p 2^{m} ( \Sigma_{j=0}^{p} (-1)^j$  ${m + j} \choose{j}$     $+\Sigma_{j=0}^{p-1} (-1)^j$ ${m + j} \choose{j}$  $)\,    -\, (-1)^p$\\

(ii) $g(p,m)=\Sigma_{j=0}^{m-1}$ ${j + p} \choose{p}$ $2^j$\\

(iii) \textit{$g(p,m)$, as a sequence parametrized by $m$, has generating function}
$$\frac{z}   {(1-z)(1-2z)^{p+1}}$$
\textit{The same expressions are valid for any polytope $P^J$ obtained from $P$ by iterated book operations, where $J$ is any $n$-tuple of non-negative integers.}\\

The last part of the theorem follows from the fact that the book construction preserves the total homology of the manifold $Z$. In particular, in the first application $P'$ is a dual-neighborly polytope of dimension $2p+1$ the homology of $Z(P')$ is concentrated in dimensions $p$ and $p+1$ and the above formula also gives its genus. This includes many odd-dimensional dual-neighborly polytopes, including all the cyclic ones.\\

Probably this is true also for any dual-neighborly polytope of odd dimension greater than 3, but we do not know yet how the genus of $Z(P)$ depends only on the combinatorics of $P$.\\

For a $3$-dimensional simple polytope $P$  the genus of $Z(P)$ is not determined by the number of facets (the simplest examples are the cube and the pentagonal book), but this is a particularity of this dimension where every simple polyhedron is dual-neighborly.\\

It must be mentioned that the experts in the field consider that, viewed from different angles, a large proportion of the simple polytopes are dual-neighborly. See [Gr], pp.129, 129a, 129b, [Z], section 4 and [P], section 1.  And these are only the roots of infinite lattices of polytopes $P^J$ stemming from them, latices that can be shown to be  disjoint for two non-combinatorially equivalent roots $P$.\\

The formula for the genus of other polytopes obtained from $P$ and its derivates by applying other operations, such as the truncation of vertices or edges (see [Gi-LdM]), can be derived from the same formulas. So the result applies to a large number of simple polytopes. We illustrate this with the case of the operation of truncating a vertex:\\

If $P$ with associated manifold $Z(P)$ is any simple $d$-polytope with $n$ facets, and $P^\vee$ the result of truncating one of its vertices, then ([Gi-LdM])
$$Z(P^\vee)=2 Z(P)	\# (2^{n-d}-1) (S^1\times S^{d-1})$$
so the genus is duplicated and then increased by $2^{n-d}-1$.\\

A simple induction shows that if the number of vertex truncations is $t$ then the genus of the resulting $Z$ is 
$$2^t (g(Z(P))-1)+t 2^{n+t-d-1}+1$$
or, in terms of the parameter $m=n-d-1$ used in section 2:
$$2^t (g(Z(P))-1)+t 2^{m+t}+1$$
showing that the result depends only on the parameters $m$ and $t$.\\

Observe that $P'$ is a $(d+1)$-polytope with $n+1$ facets, so truncating one of its vertices duplicates the genus and increases it by $2^{n-1-(d-1)}-1$. So it produces the same effect on the genus as the truncation of a vertex of $Z(P)$. Since the book construction $P'$ does not affect the genus of $Z(P)$ it follows that the effect of compositions of both operations in any order has the same effect on the genus  of $Z(P)$ as the application of the truncations only. The combinatorial type of the result of such a chain of operations on $P$ depends on the choice of facets for the book constructions, the choice of vertices to be truncated and the order of their application. The dimension of the associated manifold depends only on the number of book constructions while its genus depends only on the number of vertex truncations.\\

For example, when $P$ is a $d$-simplex, the result of $t$ vertex truncations and any number of $P'$ operations in any order is 
$$(2^t(t-1)+1)(S^1\times S^{d-1})$$ which gives the same sequence of genera $g(1,m)$ of surfaces (which can all be obtained by vertex truncation), as well as that of many $3$-dimensional polytopes for which $Z(P)$ is connected sums of copies of $S^1\times S^2$ (maybe all of them). The fact that the result of vertex truncation does not depend on the vertex chosen was first observed (in the context of intersection of real quadrics in complex space) in [B-M].\\

The effect of deeper truncations is more complicated: the genus of the result depends on the choice of simplex removed and depends on the order of the operations when combined with the book construction. Nevertheless, it can be described in terms of the combinatorics of $P$ with adequate hypotheses. We will only say, for the moment, that for many more polytopes one can obtain the genera of the associated varieties.

\section{On the sequences of genera.}

Formula (i) for the sequence of genera in Theorem 2 shows that $g(p,m)$ is of the form $2^m$ times a polinomial $R_p (m)$ of degree $p$ minus $(-1)^p$. For any value of $p$ this polynomial can be computed easily. Formulas (ii) and (iii) are more compact, but do not immediately reflect this property.\\

For $p=1$ formula (i) gives $g(1,m)=2^m(m-1)+1$, which in terms of the number of sides $n$ of the polygon ($n=m+3$) gives the usual formula $2^{n-3}(n-4)+1$. Formula (ii) gives the sum $g(1,m)=\Sigma_{j=0}^{m-1}(j + 1) 2^j$. Both formulas and the generating function appear in  https //oeis.org/A000337 together with a long list of appearances of this sequence in questions of Topology, Combinatorics, Polytope Theory and Algebraic Geometry.\\ 

For $p=2,3,4$ and $5$ the corresponding sequences also appear as entries /A055580, /A027608, /A211386 and /A21138, with their generating sequences, some formulas and various appearances in combinatorial problems.\\

Sequences for higher values of $p$ do not seem to appear. It is a curious fact that the sequences of the corresponding Euler characteristics or other variants do not seem to appear at all in the \textit{Sloane Encyclopedia of Sequences}.\\

The generating functions given in the \textit{Sloane Encyclopedia of Sequences} for those few cases gave us the clue to solve our problem. Our debt to the \textit{Encyclopedia} is partially covered by giving new formulas and a new topological interpretation to some of its sequences, as well as an infinite family of sequences generalizing them and, hopefully, by suggesting generalizations of their interpretations.\\

As a by-product of our computations we have also obtained some combinatorial identities. For example, from the formula for the number of faces of a given dimension of a neighborly polytope in [Gr] p.166, after dualizing and taking their alternating sum we get the following formula for the Euler characteristic of the manifold associated to the dual-neighborly one:

\begin{center}
$\Sigma_{k=0}^{2p} (-1)^k 2^{m+k+1} \Sigma_{j=0}^p \frac{m+2p+1}{m+k+1}$ ${m+2p-j} \choose{2p-k-j}$ ${m+k+1} \choose{2j-2p +k}$\\
\end{center}

So this expression is equal to any of the three versions of $\chi(p,m)$ given in Theorem 1 above. Perhaps this is easy to see for the experts on combinatorial identities, which is not the case of this author.\\ 

The article [O-S], which solved a topological problem by daring to deal with complicated combinatorial identities, was very stimulating for not giving up in our struggle with them.

\section*{References.}

\noindent [B-B-C-G], Bahri, A., Bendersky, M., Cohen, F.R., and Gitler, S. \emph{The polyhedral product functor: a method of computation for moment-angle complexes, arrangements and related spaces}, Adv. in Math. 22 (2010), 1634-1668.\\

\noindent [B-B-C-G2], Bahri, A., Bendersky, M., Cohen, F.R., and Gitler, S., \emph{Operations on Polyhedral Products and a new topological construction of infinite families of toric manifolds.} arXiv1011.0094v5, 2015.\\

\noindent [B-M], Bosio, F. and Meersseman, L., \textit{Real quadrics in $\C^n$, complex manifolds and convex polytopes},  Acta Mathematica, 197 (2006) 53-127.\\

\noindent [Br] Br\o ndsted, A. \textit{An introduction to convex polytopes}, 1983, Springer-Verlag.\\

\noindent [Co] Coxeter, H. (1938). Regular skew polyhedra in three and four dimension, and their topological analogues. Proceedings of the London Mathematical Society, 2(1) (1938):33-62.\\

\noindent [D-J], Davis, M. and Januszkiewicz, T., \textit{Convex polytopes, Coxeter orbifolds and torus actions}, Duke Math. Journal 62 (1991) 417-451.\\

\noindent[Gi-LdM] Gitler, S., L\'opez de Medrano, S. \textit{Intersections of Quadrics, Moment-Angle Manifolds and Connected Sums}, Geometry and Topology, 17 (2013), 1497-1534.\\

\noindent [Go], G\'omez Guti\'errez,  V., \textit{Cyclic polytopes and intersections of quadrics.}  Bolet\'in de la Sociedad Matem\'atica Mexicana, 20:2 (2014), 237?255..\\

\noindent[Go-LdM] G\'omez Guti\'errez, V., L\'opez de Medrano, S. \textit{Surfaces as Complete intersections}, in \textit{Riemann and Klein Surfaces, Automorphisms, Symmetries and Moduli Spaces},  Contemporary Math. 629 (2014), 171-180.\\

\noindent[Go-LdM2] G\'omez Guti\'errez, V., L\'opez de Medrano, S. \textit{Topology of the intersections of ellipsoids in $R^n$}, RACSAM (2018) 112: 879-891. https://doi.org/10.1007/s13398-018-0552-6,\\

\noindent [Gr], Gr\"unbaum, B., \textit{Convex polytopes}, Graduate Texts in Mathematics (221), Springer-Verlag, 2003.\\

\noindent[Hi] Hirzebruch, F., private conversation, based on
\emph{Arrangements of lines and algebraic surfaces}, in Arithmetic and Geometry, Vol. II (= Progr. Math. 36), Boston: Birkhauser, 1983, 113-140.\\

\noindent[LdM] L\'opez de Medrano, S. \textit{The Topology of the Intersection of Quadrics in} $\R^{n}$, en \textit{Algebraic Topology} (Arcata Ca,1986), Springer-Verlag Lecture Notes in Mathematics \textbf{1370} (1989), 280-292.\\

\noindent[O-S] Ortiz-Rodr\'iguez, A. and S\'anchez-Bringas, F. \textit{On some combinatorial formulae coming from Hessian Topology}, Discrete Mathematics 342(3), 2019, 615-622.\\

\noindent[P] Padrol, A. \textit{Many Neighborly Polytopes and Oriented Matroids}, Discrete Comput. Geom. 50(4), 2013, 865?902.\\

\noindent[Sl] Sloane, N.J.A., \textit{On-Line Encyclopedia of Integer Sequences (OEIS)}\newline
https://oeis.org\\
 
\noindent[Z], Ziegler, G.M., \textit{Recent Progress on Polytopes}, in \textit{Advances in Discrete and Computational
Geometry} (B. Chazelle, J.E. Goodman, R. Pollack, eds.), Contemporary Math. 223(1999), 395-406.\\
\vspace{0.3in}

\noindent\textit{santiago@matem.unam.mx\\
Instituto de Matem\'aticas.\\
Universidad Nacional Aut\'onoma de M\'exico.}\\

\noindent This work was partially supported by grant PAPIIT IN102918, UNAM.

\end{document}